\documentclass[a4paper,11pt,oneside]{article}
\usepackage[T1]{fontenc}
\usepackage[utf8]{inputenc}
\usepackage[margin=1in]{geometry}
\usepackage{amsmath,amssymb,amsthm}
\usepackage{txfonts,mathrsfs,microtype}
\usepackage[hidelinks,bookmarksnumbered]{hyperref}
\usepackage{url}
\newtheorem{theo*}{Theorem}
\newtheorem{lem}{Lemma}[section]
\newtheorem{prop}[lem]{Proposition}

\newtheorem{conj}[lem]{Conjecture}
\newtheorem*{volumerestated}{Theorem~\ref{thm:volume}}
\newtheorem*{chowpreview}{Proposition~\ref{prop:chow-lift}}
\newtheorem{rem}[lem]{Remark}
\theoremstyle{definition}
\newtheorem{defn}[lem]{Definition}
\theoremstyle{plain}

\newcommand{\C}{\mathbb C}

\newcommand{\OO}{\mathcal O}

\title{\textsc{Derived invariance of holomorphic two-forms and homological motives of stably holomorphically symplectic varieties}}
\author{\textsc{Roland Abuaf}\thanks{Rectorat de Paris, 47 rue des \'Ecoles, 75005 Paris. \textit{rabuaf@gmail.com}}}
\date{}
\begin{document}
\maketitle

\begin{abstract}
Let $\Phi:D^b(X)\to D^b(Y)$ be a derived equivalence between smooth
projective complex varieties, and let $h=c_1(L)$ for an ample line bundle
$L$ on $X$. For $z\in\mathbb C$, transporting the cup-product action and evaluating
at $e^{zh}$ gives a homomorphism of graded algebras
$\Theta^\bullet_{F,z}:H^\bullet(Y,\mathcal O_Y)\to H^\bullet(X,\mathcal O_X)$,
where $F$ is the cohomological transform of $\Phi$. We prove that
$\Theta^\bullet_{F,i}$ is an isomorphism in degrees one and two. In particular, $X$ and $Y$
have the same number of holomorphic two-forms. The proof in degree two
uses the Hodge-Riemann bilinear relations on $X$ and on $X\times X$.

We also prove that being stably holomorphically symplectic is a
derived invariant, where a smooth projective variety is called stably
holomorphically symplectic if it is a direct factor of a smooth
projective holomorphic symplectic variety. For derived-equivalent
varieties in this class, we construct an isomorphism of rational
homological motives as unital algebra objects. In particular, their
graded rational cohomology algebras are isomorphic compatibly with
the Hodge structures. Under the composition-nilpotence conjecture
for homologically trivial correspondences, the underlying isomorphism
of motives lifts to an isomorphism of rational Chow motives.
\end{abstract}

\section{Introduction}

Let $X$ and $Y$ be smooth projective complex varieties with equivalent bounded derived categories of coherent sheaves. A conjecture generally attributed to Kontsevich predicts that $X$ and $Y$ have the same Hodge numbers. In particular, they should have the same number of holomorphic two-forms. 

The derived invariance of Hochschild homology implies the invariance of the sums $\sum_{p-q=k}h^{p,q}(X)$ for all $k$. In dimension at most two, these sums, Hodge symmetry, and Poincar\'e duality determine all the Hodge numbers. In dimension three, the missing information is supplied by Popa and Schnell's theorem that the Picard varieties of derived-equivalent varieties are isogenous \cite{PS11}. In particular, $h^{1,0}$ is a derived invariant in every dimension. \\

Using generic vanishing and the Albanese morphism, Caucci and Pareschi proved the invariance of all $h^{0,q}$ for varieties of maximal Albanese dimension \cite{CP19}. For related results on cohomological support loci and the relative canonical ring of the Albanese morphism, see \cite{LP15,CLP23}.\\

In \cite{Abu17}, we proved that the graded
algebra $H^\bullet(X,\mathcal O_X)$ is a derived invariant in dimension
at most four. In arbitrary dimension, the same conclusion holds if
the Fourier-Mukai kernel is a vector bundle up to shift away from a
closed subset whose projections to both factors are proper subsets,
or if the diagonal Hodge algebras of $X$ and $Y$ are both generated
by their $(1,1)$-classes.

In dimension four, the result of \cite{Abu17} leaves only $h^{1,1}$ and $h^{2,2}$ undetermined. The dimension of $HH_0$ gives $2h^{1,1}+h^{2,2}$ once the other Hodge numbers are known. The topological signature gives $h^{2,2}-2h^{1,1}$, and its derived invariance completes the fourfold case \cite{Ant26}. In \cite{AbuSig}, we refine the signature calculations of \cite{Ant26} on each Hochschild diagonal $V_k(X)=\bigoplus_{p-q=k}H^{p,q}(X)$ and prove that $h^{n-1,1}$, $h^{n,2}$ and $h^{n-2,0}$ are derived invariants. If $K_X$ is trivial, compatibility of the Serre functors gives $K_Y\simeq\mathcal O_Y$. Serre duality then identifies $h^{0,2}$ with $h^{0,n-2}$ and gives the invariance of $h^{2,0}$ in this case. \\

Without any assumption on the canonical bundle, $H^{0,2}(X)$ lies on the same
Hochschild diagonal as $H^{k,k+2}(X)$ for $1\leq k\leq n-2$, where $n=\dim X$. To distinguish
it from these other summands, we study the graded algebra
$H^\bullet(X,\mathcal O_X)$ in more detail. Fix an equivalence $\Phi:D^b(X)\to D^b(Y)$ and an ample line bundle
$L$ on $X$. For each positive integer $m$, tensoring a class $\beta\in H^k(Y,\mathcal O_Y)$ with $\Phi(L^m)$ gives
a morphism $\beta_{\Phi(L^m)} : \Phi(L^m)\to \Phi(L^m)[k]$. Transporting this morphism
back to $X$ gives an endomorphism of degree $k$ of $L^m$. Since
$L^m$ is a line bundle, such an endomorphism is multiplication by a
unique class in $H^k(X,\mathcal O_X)$. The map used in \cite{Abu17}
is the composite
\[\Psi_m : 
 H^k(Y,\mathcal O_Y)
 \xrightarrow{\ \beta\longmapsto\beta_{\Phi(L^m)}\ }
 \operatorname{Ext}^k_Y(\Phi(L^m),\Phi(L^m))
 \xrightarrow{ \Phi^{-1} }
 \operatorname{Ext}^k_X(L^m,L^m)
 \simeq H^k(X,\mathcal O_X).
\]
This map respects products because scalar multiplication and transport by
$\Phi^{-1}$ respect composition. Let us write this composite in cohomological terms. Denote by $F$ the
cohomological Fourier-Mukai transform of $\Phi$, put $h=c_1(L)$ and
$s_X=\sqrt{\operatorname{td}(X)}$, and recall that the Mukai vector
of $L^m$ is $v(L^m)=s_Xe^{mh}$. The functor $\Phi$ takes
multiplication by $\Psi_m(\beta)$ on $L^m$ to multiplication by $\beta$ on
$\Phi(L^m)$. By compatibility of the Hochschild action with the Chern
character \cite{Cal05} and the compatibility with evaluation on
objects \cite[Corollary~6.3.2]{BF08}, this gives
\[
             F\bigl(\Psi_m(\beta) \cup v(L^m)\bigr)=\beta\cup F\bigl(v(L^m)\bigr).
\]
Since $v(L^m)$ has constant term one, it is invertible in the
cohomology ring. We therefore obtain
\[
 \Psi_m(\beta) =(s_Xe^{mh})^{-1}\cup
        F^{-1}\!\left(\beta\cup F(s_Xe^{mh})\right).
\]

\textbf{Our main goal is to prove that, under suitable assumptions on $k$, $m$, and $h$, a variant of the map $\beta \longmapsto \Psi_m(\beta)$ is injective.} \\

The strategy in \cite{Abu17} was to find a line bundle whose image
under $\Phi$ has nonzero rank. For the family considered here, it
would suffice to find a positive integer $m$ for which $\Phi(L^m)$
has nonzero rank. Indeed, if $\Psi_m(\beta)=0$, then
$\beta_{\Phi(L^m)}=0$, and taking the trace gives
$\operatorname{rk}(\Phi(L^m))\beta=0$. Such an $m$ therefore gives
injectivity in every degree. The difficulty is to establish the existence of such a line bundle
for a general equivalence. \\

We modify the map $\Psi_m$ by replacing the Mukai vector $s_Xe^{mh}$ of
$L^m$ by its Chern character $e^{mh}$ and allowing the parameter $m$ to take
complex values. For $z\in\mathbb C$ and $\beta\in H^{\bullet} (Y,\mathcal O_Y)$, put
\begin{equation}\label{eq:evaluation-map}
 \begin{aligned}
 \Theta^\bullet_{F,z}:H^\bullet(Y,\mathcal O_Y)&\longrightarrow
                                  H^\bullet(X,\mathcal O_X),\\
 \beta&\longmapsto e^{-zh}\cup F^{-1}\!\left(\beta\cup F(e^{zh})\right).
 \end{aligned}
\end{equation}

The main reason for replacing the Mukai vector of $L^m$ by its Chern
character is that it greatly simplifies the cohomological calculations (see
Section~\ref{subsec:positivity}). \textbf{We wish to find a value of $z$ for which $\Theta^2_{F,z}$ is injective.} Let us first explain the strategy in degree one. Replacing the categorically defined map $\Psi_m : H^{\bullet}(Y, \OO_Y) \longrightarrow H^{\bullet}(X, \OO_X)$ by its cohomological variant $\Theta^{\bullet}_{F,z}$ and allowing $z$ to vary in $\C$ enables us to exploit orthogonality relations between the terms contributing to $\Theta^{1}_{F,z}$.



More precisely, assume $n\geq1$ and let $\beta\in H^{0,1}(Y)$ and $D_\beta=F^{-1}(\beta \cup -)F : H^{\bullet}(X, \C) \longrightarrow H^{\bullet}(X, \C)$. The compatibility of the Todd-corrected Hochschild-Kostant-Rosenberg isomorphisms with the module action (see \cite{CRV12}) allows us to write $D_{\beta} = m_a + \iota_{v}$, where $a \in H^{1}(X,\OO_X)$, $v \in H^0(X, T_X)$, $m_a$ is multiplication by $a$, and $\iota_v$ is contraction with the vector field $v$. Graded commutativity gives $\beta \cup \overline{\beta}+\overline{\beta}\cup\beta=0$. Since $F$ commutes with complex conjugation, we obtain
\[D_\beta\overline{D_\beta}+\overline{D_\beta}D_\beta=0,\] where
$\overline{D_\beta}(x)=\overline{D_\beta(\bar x)}$. Evaluating this identity at $1_X$ gives $\iota_{v}(\overline{a}) + \overline{\iota_{v}(\overline{a})} = 0$, so $\iota_v(\overline{a})$ is purely imaginary. On the other hand, the Hodge-Riemann bilinear relations in degree $1$ imply that the Hermitian form
\[\langle u,w \rangle = i \int_{X} \overline{u} \cup w \cup h^{n-1}\]
is positive definite on $H^{0,1}(X)$. We have $\overline{a} \cup h^n = 0$ for reasons of Hodge type. Applying $\iota_{v}$ and using the graded Leibniz rule gives
\[\iota_v(\overline{a}) \cup h^n = n \,\overline{a}\cup\iota_v(h)\cup h^{n-1}.\]
Integrating both sides and using the fact that $\iota_v(\overline{a})$ is purely imaginary, we obtain
\[ \operatorname{Im} \langle a , \iota_v(h) \rangle = 0,\]
and hence $\operatorname{Re} \langle a, i \iota_v(h) \rangle = 0$. Since $D_{\beta} = m_a+ \iota_v$, a direct computation gives
\[ \forall z \in \C, \ \Theta^1_{F,z}(\beta) = a + z \iota_v(h).\]
Thus $a$ and $i\iota_v(h)$ are orthogonal for the real part of the Hermitian product $\langle \cdot, \cdot \rangle$, and
\[ \lVert \Theta^1_{F,i}(\beta) \rVert^2 = \lVert a \rVert^2 + \lVert \iota_v(h) \rVert^2,\]
where $\lVert \cdot \rVert$ is the associated norm. If $\Theta^1_{F,i}(\beta)=0$, this equality gives $a=0$ and $\iota_v(h)=0$. A standard Lefschetz argument then gives $D_\beta=0$, hence $\beta=0$. This proves the injectivity of $\Theta^1_{F,i}$. \\

These ideas and techniques extend to the case $k=2$. The following is our first main result.

\begin{theo*}\label{thm:twoforms}
Let $\Phi:D^b(X)\to D^b(Y)$ be a derived equivalence of smooth projective
complex varieties, let $F$ be its cohomological transform, and let
$L$ be an ample line bundle on $X$. For $h=c_1(L)$, the homomorphism
$\Theta^\bullet_{F,i}$ is an isomorphism in degrees one and two. In
particular, $h^{1,0}(X)=h^{1,0}(Y)$ and $h^{2,0}(X)=h^{2,0}(Y)$.
\end{theo*}

The proof in degree two begins with the same idea. For
$\beta\in H^{0,2}(Y)$, write $D_\beta=F^{-1}m_\beta F=m_a+K+\Pi$,
where $a\in H^{0,2}(X)$, $K$ is a Kodaira-Spencer contraction,
and $\Pi$ is contraction by a bivector. Put $\eta=K(h)$ and
$c=\Pi(h^2)/2$. Evaluation at $i$ gives
$\Theta^2_{F,i}(\beta)=a+i\eta-c$. In this even degree, $D_\beta$ commutes with its
complex conjugate. As in degree one, this relation yields
orthogonality for the real part of the Hodge-Riemann product,
and hence
\[
 \|\Theta^2_{F,i}(\beta)\|_h^2
       =\|a-c\|_h^2+\|\eta\|_h^2,
\]
where the norm is now the positive Hodge-Riemann norm on
$H^{0,2}(X)$. Thus $\Theta^2_{F,i}(\beta)=0$ implies $K(h)=0$
and $a=c$. The new step is to apply positivity to the operators
themselves. We represent the actions of $K$ and $\Pi$ on even
cohomology by correspondences on $X\times X$. After a sign change
on the second factor, their commutator identities with
multiplication by $h$ become primitivity conditions for the
ample class $p_1^*h+p_2^*h$. The Hodge-Riemann relations then
give nonnegative trace expressions whose sum is
$\operatorname{Tr}(D_\beta\overline{D_\beta})$.
This trace vanishes because $D_\beta$ is nilpotent, as it lowers
the Hochschild grading by two, and commutes with its conjugate.
Strict positivity forces $K=\Pi=0$ on even cohomology, and
$a=c$ then gives $a=0$. Transporting back to $Y$ and evaluating
at $1_Y$ yields $\beta=0$. \\

The properties of the morphism $\Theta^{\bullet}_{F,i}$ lead to a second result, concerning
the cohomology algebra and its realization by algebraic correspondences.
We call a smooth projective complex variety $X$ \emph{stably
holomorphically symplectic} if there exists a smooth projective complex
variety $Z$ such that $X\times Z$ is holomorphically symplectic.
This class contains all holomorphic symplectic varieties and also
odd-dimensional varieties, such as the product of a K3 surface with
an elliptic curve. Proposition~\ref{prop:stable-derived} shows that
being stably holomorphically symplectic is a derived invariant.

For this entire class, we construct an isomorphism of rational
homological motives that preserves the unit and cup product.
We write $h_{\mathrm{hom}}(Z)_{\mathbb Q}$ for the motive of $Z$
with rational coefficients, modulo homological equivalence defined
by Betti cohomology.

\newcommand{\volumestatement}{%
Let $X$ be a smooth projective stably holomorphically symplectic
complex variety, and let $Y$ be a smooth projective variety derived
equivalent to $X$. Then $Y$ is stably holomorphically symplectic, and
\[
 \bigl(h_{\mathrm{hom}}(X)_{\mathbb Q},\cup,1_X\bigr)
       \simeq
 \bigl(h_{\mathrm{hom}}(Y)_{\mathbb Q},\cup,1_Y\bigr)
\]
as unital algebra objects. Consequently, there is an isomorphism of
graded $\mathbb Q$-algebras
$H^\bullet(X,\mathbb Q)\simeq H^\bullet(Y,\mathbb Q)$ whose
component in every degree is an isomorphism of rational Hodge
structures. In particular, $h^{p,q}(X)=h^{p,q}(Y)$ for all $p,q$.}
\begin{theo*}\label{thm:volume}
\volumestatement
\end{theo*}

The theorem gives an isomorphism that is rational, preserves
the Hodge structures and cup product, and is induced by an algebraic
correspondence with an inverse modulo homological equivalence.
It applies in every dimension, without a restriction on deformation
type or on the structure of the LLV algebra. In particular, it proves
the homological form of the multiplicative Orlov conjecture for
hyperk\"ahler varieties \cite[Conjecture~0.3]{MSY26} and extends its
conclusion to all stably holomorphically symplectic varieties.
Fu and Vial formulate a Chow-theoretic conjecture for Frobenius
algebra objects \cite[Conjecture~4.7]{FV21}.

Several parts of this conclusion were known in important cases.
For hyperk\"ahler varieties, Taelman proved that derived-equivalent
varieties have isomorphic rational Hodge structures in every degree
\cite[Theorem~D]{Tae23}; see also Beckmann's exposition
\cite[Section~8, Theorem~8.2]{Bec22}. For holomorphic symplectic varieties
beyond the hyperk\"ahler case, Leijnse established the rational
Hodge-theoretic conclusion in dimension four and, in higher dimensions,
under additional hypotheses on the LLV algebras
\cite[Theorems~1.8 and~7.21]{Lei26}. Maulik, Shen and Yin proved
the isomorphism of homological motives preserving cup product for
varieties of $K3^{[d]}$-type \cite[Theorem~0.4]{MSY26}.
Theorem~\ref{thm:volume} extends this homological cup-product
statement to all hyperk\"ahler varieties, all holomorphic symplectic
varieties, and their smooth projective direct factors.\\

We first prove the equality of Hodge numbers directly. The geometric
characterization in Proposition~\ref{prop:stable-characterization}
shows that a stably holomorphically symplectic variety is holomorphic
symplectic in even dimension and becomes so after taking the product
with an elliptic curve in odd dimension. In the even-dimensional
case, the morphism $\Theta^2_{F,i}$ transports a symplectic two-form
to a symplectic two-form on the derived partner. Comparing the ranks
of their powers on the Hochschild diagonals then recovers all Hodge
numbers. The odd-dimensional case follows by the K\"unneth formula.

The passage from Hodge numbers to motives rests on a normalization
of the cohomological transform. Suppose first that $\dim X=2d>0$,
and let $\tau\in H^{0,2}(X)$ be conjugate to a holomorphic symplectic
form. Choose $\beta\in H^{0,2}(Y)$ with
$\Theta^2_{F,i}(\beta)=\tau$. By multiplicativity,
$\Theta^{2d}_{F,i}(\beta^d)=\tau^d\ne0$. Since $F(e^{zh})$ has
diagonal Hodge type, we have
\[
 \beta^d\cup F(e^{zh})=[F(e^{zh})]_0\,\beta^d,
\]
where $[-]_j$ denotes the component of cohomological degree $j$.
Thus the rational polynomial $r(z)=[F(e^{zh})]_0$ is nonzero at
$i$, and hence at some positive integer $m$.

The derived invariance of the Looijenga-Lunts-Verbitsky algebra,
proved by Taelman, then allows us to write
\[
 F(e^{mh})=r(m)e^b,
 \qquad b\in H^2(Y,\mathbb Q)\cap H^{1,1}(Y).
\]
The normalized transform $G=m_{e^{-b}}F m_{e^{mh}}$ therefore
preserves the unit line. The transport of the degree-two Hochschild
action, together with an elementary symplectic Lefschetz identity,
then shows that $G$ and its inverse are triangular in cohomological
degree. Their degree-preserving
components are inverse and are induced by algebraic correspondences
of codimension $\dim X$. Writing $G_0$ for the degree-preserving
component, we obtain an isomorphism of homological motives whose
action $U=r(m)^{-1}G_0$ sends $1_X$ to $1_Y$.

To prove multiplicativity, we use the full Hochschild action.
The transform $G$ transports all contraction operators, and its
triangularity implies that $U$ still transports each homogeneous
contraction component. A holomorphic symplectic form $\sigma$ lets
us recover every cup-product operator from these contractions:
locally,
\[
                    [m_\sigma,\iota_v]
                        =-m_{\iota_v\sigma},
\]
and iterating this identity gives multiplication by forms of every
degree. Thus $Um_xU^{-1}$ is again a multiplication operator;
evaluation at $1_Y$ identifies it with $m_{U(x)}$. This proves
$U(x\cup y)=U(x)\cup U(y)$.
In odd dimension we apply the construction to $X\times E$ and
$Y\times E$, with $E$ an elliptic curve and a product polarization.
The resulting map factors as $U_X\otimes\operatorname{id}_E$,
so both the motivic isomorphism and its multiplicativity descend
to $X$ and $Y$.\\

The correspondence inducing this isomorphism is defined in the Chow
group of $X\times Y$. Since it is invertible modulo homological
equivalence, it is natural to ask whether it already admits an inverse
modulo rational equivalence. An affirmative answer would establish,
for the varieties considered here, Orlov's conjecture that
derived-equivalent smooth projective varieties have isomorphic
rational Chow motives \cite[Conjecture~1]{Orl05}.
Write $h_{\mathrm{CH}}(Z)_{\mathbb Q}$ for the motive of $Z$ with
rational coefficients modulo rational equivalence. For the varieties
in Theorem~\ref{thm:volume}, Orlov's conjecture follows from
composition nilpotence for homologically trivial correspondences.
As we recall in Section~\ref{subsec:chow}, this nilpotence is
implied by any one of the conjectures of Kimura-O'Sullivan on finite
dimensionality, of Voevodsky on smash nilpotence, and of
Bloch-Beilinson on the filtration of Chow groups.

\newcommand{\chowliftstatement}{%
Let $X$ and $Y$ be derived-equivalent smooth projective complex
varieties, with $X$ stably holomorphically symplectic. If the composition-nilpotence conjecture
for homologically trivial correspondences holds for $X$ and $Y$ (see Conjecture~\ref{conj:composition-nilpotence}), the underlying isomorphism of homological motives in
Theorem~\ref{thm:volume} lifts to an isomorphism of rational Chow motives
\[
       h_{\mathrm{CH}}(X)_{\mathbb Q}\simeq h_{\mathrm{CH}}(Y)_{\mathbb Q}.
\]
In particular, Orlov's conjecture holds for the pair $(X,Y)$.}
\begin{chowpreview}
\chowliftstatement
\end{chowpreview}

Beyond these motivic applications, the methods developed here reveal
a rich interplay between Hochschild cohomology and Hodge-Riemann
positivity. Relations preserved by derived equivalences acquire
geometric content through orthogonality and positivity, both on a
variety and on its products. The range of consequences obtained
from this interaction suggests that it may provide a fruitful
approach to further questions of derived invariance. In particular,
the graded algebra structure of our construction invites us to
explore how much of this method extends beyond degree two.

Degree three presents the first substantial challenge. The trace
calculation now involves mixed terms whose signs are harder to
control, and the choice of polarization becomes significant, as
illustrated by Remark~\ref{rem:degree-three-polarization}. These
difficulties point towards a closer study of the positivity
relations among Hochschild operators and their dependence on the
polarization. We believe that the flexibility of this approach,
together with the wealth of relations supplied by the Hochschild
action, offers promising possibilities for progress on the derived
invariance of $h^{0,3}$ and, more broadly, of Hodge numbers in higher
degrees.\\

In Section~\ref{subsec:graded}, we construct the homomorphisms
$\Theta^\bullet_{F,z}$ and prove that $\Theta^1_{F,i}$ is injective. Sections~\ref{subsec:positivity}
and~\ref{subsec:evaluation} establish the positivity statements and
use them to prove injectivity in degree two. In
Section~\ref{subsec:hodge-numbers} we characterize stably holomorphically
symplectic varieties, prove that this property is derived invariant,
and compare symplectic Lefschetz ranks to recover all Hodge numbers.
Section~\ref{subsec:motives} constructs the isomorphism of homological
motives and proves its compatibility with cup product, completing
the proof of Theorem~\ref{thm:volume}. Section~\ref{subsec:chow}
gives the conditional lift to Chow motives.

\bigskip

\noindent \textbf{Acknowledgements and AI disclosure.} During the preparation of this article, the author used ChatGPT (OpenAI) as an interactive tool for mathematical discussion, verification of computations, exploration of proof strategies, assistance in locating relevant references, and improvement of the English. The paper was written by the author, who is responsible for its correctness. \\

\section{The graded homomorphism and derived invariance of two-forms}\label{sec:twoforms}

We prove Theorem~\ref{thm:twoforms} in this section. We first
construct the graded algebra homomorphism $\Theta^\bullet_{F,z}$
and treat the degree-one case, which serves as a model for the proof
in degree two. The rest of the section is devoted to this latter
case. In Section~\ref{subsec:positivity} we establish the
orthogonality and positivity results we need, and in
Section~\ref{subsec:evaluation} we use them to prove the injectivity
of $\Theta^2_{F,i}$.

\subsection{Construction of the homomorphism and injectivity in degree one}\label{subsec:graded}

Let $\Phi:D^b(X)\xrightarrow{\sim}D^b(Y)$ be an equivalence. We first
recall how to transport cup product by a class in $H^k(Y,\mathcal O_Y)$
to an operator on the full cohomology of $X$. The use of the full cohomology
is essential here, since for odd $k$ this operator exchanges the even and odd parts.

Put
\[
 HT^j(X)=\bigoplus_{r+s=j}H^s(X,\wedge^r T_X)
 \quad \textrm{and} \quad s_X=\sqrt{\operatorname{td}(X)},
\]
where the square root has constant term one. We use the Todd-corrected
Hochschild-Kostant-Rosenberg isomorphisms
\[
 \begin{aligned}
 I_X:HH_k(X)&\xrightarrow{\ \sim\ }
       \bigoplus_{p-q=k}H^q(X,\Omega_X^p),\\
 J_X:HT^j(X)&\xrightarrow{\ \sim\ }HH^j(X).
 \end{aligned}
\]
On homology, $I_X$ is the usual HKR map followed by cup product with
$s_X$. On Hochschild cohomology, $J_X$ is obtained by first contracting polyvectors
with $s_X$ and then applying the cohomological HKR map. 

For a class $\gamma\in H^s(X,\wedge^r T_X)$, let
\[
 \iota_\gamma:H^q(X,\Omega_X^p)\longrightarrow
                    H^{q+s}(X,\Omega_X^{p-r})
\]
be the operator obtained by taking the cup product with $\gamma$ in
sheaf cohomology and applying the contraction
$\wedge^r T_X\otimes\Omega_X^p\to\Omega_X^{p-r}$. This operator vanishes
for $p<r$. Our convention on decomposable polyvectors is
$\iota_{v_1\wedge\cdots\wedge v_r}=\iota_{v_r}\cdots\iota_{v_1}$.
Cup products and contractions use the usual Koszul signs for
total cohomological degree.
In particular, when $r=0$, the operator $\iota_\gamma$ is ordinary cup
product by $\gamma\in H^s(X,\mathcal O_X)$ and we denote it by $m_{\gamma}$. The compatibility theorem
of Calaque, Rossi and Van den Bergh \cite[Theorem~1.4]{CRV12} identifies
this action with the Hochschild module action. Namely, for
$w\in HH_\bullet(X)$,
\[
             I_X\bigl(J_X(\gamma)\cap w\bigr)=\iota_\gamma I_X(w).
\]

Let $\mathcal P$ be a Fourier-Mukai kernel representing $\Phi$
\cite{Huy06}, and let $p_X,p_Y$ be the projections from $X\times Y$.
By \cite[Theorem~1.2]{MS09}, the Todd-corrected HKR maps $I_X$ and
$I_Y$ identify the action $\Phi_*$ on Hochschild homology with the
Fourier-Mukai transform $F=I_Y\Phi_*I_X^{-1}$ on total cohomology,
given by
\begin{equation}\label{eq:mukai}
 F(x)=p_{Y*}\!\left(
       p_X^*x\cup\operatorname{ch}(\mathcal P)
       \cup p_X^*s_X\cup p_Y^*s_Y\right).
\end{equation}
We use the Betti normalizations of characteristic classes and Gysin maps.
The characteristic classes in this formula are rational and of diagonal
Hodge type. The Gysin map is real and lowers
both Hodge indices by $\dim X$. Consequently $F$ commutes with complex
conjugation and preserves the difference $p-q$. Since $p+q$ and $p-q$
have the same parity, $F$ also identifies
$V_X=H^{\mathrm{ev}}(X,\mathbb C)$ with $V_Y=H^{\mathrm{ev}}(Y,\mathbb C)$. \\

It remains to describe how $F$ acts on the operators $\iota_\gamma$.
Write $\Phi^{HH}:HH^\bullet(X)\to HH^\bullet(Y)$ for the isomorphism
on Hochschild cohomology. Given $\gamma_X\in HT^\bullet(X)$, let
$\gamma_Y=J_Y^{-1}\bigl(\Phi^{HH}(J_X(\gamma_X))\bigr)$. Functoriality
of the Hochschild module action gives
\[
 \Phi_*\bigl(J_X(\gamma_X)\cap w\bigr)
       =J_Y(\gamma_Y)\cap\Phi_*w.
\]
Applying $I_Y$ and using the compatibility with contraction on both
sides gives the transport identity
\begin{equation}\label{eq:transport}
 F\,\iota_{\gamma_X}F^{-1}=\iota_{\gamma_Y}
             \quad\text{on }H^\bullet(Y,\mathbb C).
\end{equation}
In particular, for $\beta\in H^k(Y,\mathcal O_Y)$, the operator
$D_\beta=F^{-1}m_\beta F$ is the action of a class in $HH^k(X)$.

Fix an ample line bundle $L$ on $X$ and set $h=c_1(L)$. For
$z\in\mathbb C$ and $\beta\in H^k(Y,\mathcal O_Y)$, write
\[
 T_{\beta,z}=m_{e^{-zh}}D_\beta m_{e^{zh}}.
\]

\begin{prop}\label{prop:graded-homomorphism}
For every $z\in\mathbb C$, the formula
\[
 \begin{aligned}
 \Theta^\bullet_{F,z}:H^\bullet(Y,\mathcal O_Y)&\longrightarrow
                                      H^\bullet(X,\mathcal O_X),\\
 \beta&\longmapsto T_{\beta,z}(1_X)
   =e^{-zh}\cup F^{-1}\!\left(\beta\cup F(e^{zh})\right)
 \end{aligned}
\]
defines a unital homomorphism of graded $\mathbb C$-algebras.
Moreover, for all $k,j\geq0$, $\beta\in H^k(Y,\mathcal O_Y)$ and
$u\in H^j(X,\mathcal O_X)$, one has
\begin{equation}\label{eq:first-column}
 T_{\beta,z}(u)=\Theta^k_{F,z}(\beta)\cup u
                    \in H^{k+j}(X,\mathcal O_X).
\end{equation}
\end{prop}

\begin{proof}
Let
\[
 j_Y:H^\bullet(Y,\mathcal O_Y)\hookrightarrow HT^\bullet(Y)
 \quad \textrm{and} \quad
 \operatorname{pr}_0:HT^\bullet(X)\longrightarrow
                         H^\bullet(X,\mathcal O_X)
\]
denote the inclusion and projection corresponding to polyvector
degree zero. Both are unital homomorphisms of graded algebras:
for the projection, this follows because the components of positive
polyvector degree form an ideal. By \cite[Theorem~1.4]{CRV12},
the corrected HKR maps $J_X$ and $J_Y$ are also isomorphisms of
graded algebras.

For an integer $m$, consider the equivalence
$\Phi_m=\Phi\circ(-\otimes L^m)$, whose cohomological transform is
$F_m=F\circ m_{e^{mh}}$. Thus
\[
 T_{\beta,m}=F_m^{-1}m_\beta F_m.
\]
By the transport identity \eqref{eq:transport}, this operator is
contraction by
\[
 \gamma_{\beta,m}
   =J_X^{-1}(\Phi_m^{HH})^{-1}J_Yj_Y(\beta).
\]
Every component of positive polyvector degree vanishes on
$H^\bullet(X,\mathcal O_X)$, while the component of polyvector
degree zero acts by cup product. Evaluating at $1_X$ therefore gives
\[
 \Theta^\bullet_{F,m}
   =\operatorname{pr}_0\circ J_X^{-1}
      \circ(\Phi_m^{HH})^{-1}\circ J_Y\circ j_Y.
\]
Every map in this composition is a unital homomorphism of graded
algebras, so the same holds for $\Theta^\bullet_{F,m}$.
The same description gives, for
$u\in H^j(X,\mathcal O_X)$,
\[
 T_{\beta,m}(u)
   =\operatorname{pr}_0(\gamma_{\beta,m})\cup u
   =\Theta^k_{F,m}(\beta)\cup u,
\]
which proves \eqref{eq:first-column} for every integer $m$.

Finally, since $h$ is nilpotent, $T_{\beta,z}$ and
$\Theta^k_{F,z}(\beta)$ depend polynomially on $z$.
For fixed homogeneous inputs, preservation of degree, products
and the unit, as well as \eqref{eq:first-column}, are polynomial
identities. To express the degree condition as a polynomial identity, we pass to
$H^\bullet(X,\mathbb C)/H^k(X,\mathcal O_X)$.
These identities hold at every integer, hence for every
$z\in\mathbb C$. Complex linearity is immediate from the definition.
\end{proof}

As a warm-up for the proof of injectivity in degree two, we first
treat degree one. This simpler case already illustrates how to use
the relation between a transported operator and its complex conjugate
to bring the Hodge-Riemann bilinear relations into play.

\begin{prop}\label{prop:degree-one}
For every ample line bundle $L$ on $X$ and $h=c_1(L)$, the map
$\Theta^1_{F,i}:H^1(Y,\mathcal O_Y)\to H^1(X,\mathcal O_X)$ is an
isomorphism. In particular, $h^{1,0}(X)=h^{1,0}(Y)$.
\end{prop}
\begin{proof}
Let $\beta\in H^1(Y,\mathcal O_Y)$. By \eqref{eq:transport} and the
decomposition of $HT^1(X)$, its transported action has the form
\[
 D_\beta=m_a+\iota_v,\quad \textrm{with} \
 a\in H^{0,1}(X) \ \textrm{and} \ v\in H^0(X,T_X).
\]
If $\dim X=0$, both summands vanish, so $D_\beta=0$ and hence
$\beta=0$. Injectivity is therefore immediate in this case.
Assume now that $n=\dim X\geq1$, and put
$b=\iota_v(h)\in H^{0,1}(X)$.
Since $\iota_v$ is an odd derivation of total degree $-1$, we have
$\iota_v(e^{zh})=z e^{zh}b$, and therefore
\[
                    \Theta^1_{F,z}(\beta)=a+zb.
\]
The classes $\beta$ and $\bar\beta$ have odd degree, so their
multiplication operators anticommute. Since $F$ commutes with complex
conjugation, their transported operators also anticommute:
\[
 D_\beta\overline{D_\beta}+\overline{D_\beta}D_\beta=0 \]
 where $\overline{D_\beta}(x)=\overline{D_\beta(\bar x)}$. Applying this identity to $1_X$ and using
$a\bar a+\bar a a=0$, we obtain
\[
 c+\bar c=0,\quad \textrm{where} \
 c=\iota_v(\bar a)\in H^{0,0}(X)=\mathbb C.
\]
Thus $c$ is purely imaginary. On the other hand, $\bar a h^n=0$
by type. The graded Leibniz rule gives
\[
 0=\iota_v(\bar a h^n)
   =c h^n-n\bar a\,b\,h^{n-1}.
\]
Consider the positive definite Hermitian form on $H^{0,1}(X)$
\[
 \langle u,w\rangle_h
       =i\int_X\bar u\,w\,h^{n-1},
 \quad \textrm{with} \ \|u\|_h^2=\langle u,u\rangle_h,
\]
whose positivity follows from the degree-one Hodge-Riemann relations
\cite[p.~123]{GH78}. Integrating the preceding identity shows that
\[
 \langle a,b\rangle_h
       =\frac{ic}{n}\int_X h^n\in\mathbb R.
\]
Consequently $\operatorname{Re}\langle a,ib\rangle_h=0$. Expanding
the square of the norm therefore gives
\begin{equation}\label{eq:degree-one-pythagoras}
 \|\Theta^1_{F,i}(\beta)\|_h^2
      =\|a+ib\|_h^2
      =\|a\|_h^2+\|b\|_h^2.
\end{equation}

Suppose now that $\Theta^1_{F,i}(\beta)=0$. By
\eqref{eq:degree-one-pythagoras}, we have $a=b=0$.
It remains to show that $\iota_v$ vanishes on cohomology.
Since $\iota_v(h)=b=0$, this derivation commutes with $L_h=m_h$.
If $x$ is primitive of degree $1\leq j\leq n$, then
\[
 L_h^{n-j+1}\iota_v(x)
       =\iota_v(L_h^{n-j+1}x)=0.
\]
Hard Lefschetz identifies $L_h^{n-j+1}$ with an isomorphism from
$H^{j-1}(X,\mathbb C)$ to $H^{2n-j+1}(X,\mathbb C)$, so
$\iota_v(x)=0$. The operator also vanishes on $H^0(X,\mathbb C)$.
Since every cohomology class is a sum of Lefschetz powers of
primitive classes, commutation with $L_h$ now gives $\iota_v=0$.

Thus $D_\beta=0$. Conjugating by $F$ gives $m_\beta=0$, and
evaluation at $1_Y$ yields $\beta=0$. This proves injectivity.
Applying the same argument to the inverse equivalence, with an ample
line bundle on $Y$, gives the reverse dimension inequality.
Hence $\Theta^1_{F,i}$ is an isomorphism, and Hodge symmetry gives
$h^{1,0}(X)=h^{1,0}(Y)$.
\end{proof}

\subsection{Contraction identities and Hodge-Riemann positivity}\label{subsec:positivity}

We turn to the degree-two component of $\Theta^\bullet_{F,i}$.
As in degree one, we first use the relation between a transported
operator and its conjugate to obtain orthogonality for the real part
of a positive Hermitian form. The resulting norm identity will show
that $K(h)=0$ whenever $\Theta^2_{F,i}(\beta)=0$. We then use
Hodge-Riemann positivity for the correspondences of $K$ and $\Pi$
on $X\times X$ to prove that both operators vanish. \\

Denote
by $\mathcal V=H^\bullet(X,\mathbb C)$ the full Hodge cohomology of $X$,
by $V=H^{\mathrm{ev}}(X,\mathbb C)$ its even part, and by
\[
 \mathscr E_X=\operatorname{im}\!\left(
 HH^2(X)\longrightarrow\operatorname{End}_{\mathbb C}(V)\right)
\]
the image of the degree-two Hochschild action on $V$. In this
subsection, assume $n=\dim X\geq2$ and fix $h=c_1(L)$ for an ample
line bundle $L$ on $X$. For the operator calculations, write $A_u=m_u$
and $L_h=A_h$. On $\mathcal V$ a degree-two Hochschild class acts as
\[
 D=A_a+K+\Pi,\quad \textrm{with} \  K=\iota_\kappa \ \textrm{and} \ \Pi=\iota_\pi,
\]
where $a\in H^{0,2}(X)$, $\kappa\in H^1(X,T_X)$ and $\pi\in H^0(X,\wedge^2T_X)$. The three summands have bidegrees $(0,2)$, $(-1,1)$ and $(-2,0)$, respectively. The restrictions of these three components to $V$ are determined by $D|_V$, since they have distinct shifts in total degree. We use the bilinear Poincar\'e pairing $P(x,y)=\int_X xy$ on $V$ and write $T^t$ for the transpose of $T$ with respect to $P$. The conjugate operator is $\overline T(x)=\overline{T(\overline x)}$. These two operations commute because $P$ is real.

\begin{samepage}
\begin{lem}\label{lem:contraction}
On the full Hodge cohomology $\mathcal V$ one has
\begin{equation}\label{eq:operators-L}
 [L_h,K]=-A_\eta,\qquad [L_h,[L_h,\Pi]]=2A_c,\qquad [L_h,[L_h,[L_h,\Pi]]]=0,
\end{equation}
where $ \eta=K(h),\  c=\tfrac12\Pi(h^2),\ \textrm{and} \ \eta,c\in H^{0,2}(X)$. Moreover,
\begin{equation}\label{eq:operators-exp}
 D(e^{zh})=e^{zh}(a+z\eta+z^2c)\qquad(z\in\mathbb C).
\end{equation}
On $V$, the transposes of $A_a,K,\Pi$ are $A_a,-K,\Pi$, respectively.

\end{lem}
\end{samepage}
\begin{proof}

Contraction by a class in $H^1(X,T_X)$ is a derivation of total
degree zero. Thus
\[
 [L_h,K](x)=hK(x)-K(hx)=-K(h)x,
\]
which proves $[L_h,K]=-A_\eta$. To compute $J=[L_h,\Pi]$, let $b_\pi:\Omega_X^1\to T_X$ be the
bundle map defined by
\[
 \lambda\bigl(b_\pi(\theta)\bigr)=\pi(\lambda,\theta).
\]
With our contraction convention, the local computation in the proof
of \cite[Lemma~2.13]{Tae23} gives
\[
 J=\iota_\xi,
 \quad \textrm{with} \
 \xi=H^1(b_\pi)(h)\in H^1(X,T_X).
\]
Although Taelman assumes that the bivector is nondegenerate,
this computation uses only the exterior Leibniz rule and remains
valid for an arbitrary bivector. In particular, $J$ is also a
derivation of total degree zero. Since $\Pi(h)=0$ by bidegree,
\[
 J(h)=h\Pi(h)-\Pi(h^2)=-2c.
\]
Applying the same argument to $J$ gives
\[
 [L_h,[L_h,\Pi]]=[L_h,J]=-A_{J(h)}=2A_c.
\]
Multiplication by $c$ commutes with $L_h$, so the third iterated
commutator vanishes. This proves \eqref{eq:operators-L}.\\

We now compute the conjugated operator. The identities just proved,
together with $[L_h,A_a]=0$, give
\[
 [L_h,D]=-A_\eta+J,\qquad
 [L_h,[L_h,D]]=2A_c,\qquad
 \operatorname{ad}_{L_h}^3(D)=0.
\]
Consequently the polynomial $e^{-zL_h}De^{zL_h}$ has zero third
derivative. Evaluating its first two derivatives at zero yields
\[
 \begin{aligned}
 e^{-zL_h}De^{zL_h}
   &=D-z[L_h,D]+\frac{z^2}{2}[L_h,[L_h,D]]\\
   &=A_{a+z\eta+z^2c}+(K-zJ)+\Pi.
 \end{aligned}
\]
The contractions $K,J,\Pi$ vanish on $1_X$. Applying this identity to
$1_X$ therefore gives
\[
 e^{-zh}\cup D(e^{zh})=a+z\eta+z^2c,
\]
which is \eqref{eq:operators-exp}. Applied to $D=\Pi$, the same
identity gives $\Pi(e^{zh})=z^2e^{zh}c$; comparison of coefficients
also yields
\[
 \Pi(h^m)=m(m-1)h^{m-2}c\qquad(m\geq2).
\]

It remains to determine the transposes on $V$.
Multiplication by an even class is symmetric for $P$, so
$A_a^t=A_a$ and $L_h^t=L_h$. Since $K$ is a derivation,
\[
 P(Kx,y)+P(x,Ky)=\int_X K(xy)=0.
\]
Indeed, only products of total degree $2n$ can contribute, and
$K$ vanishes on $H^{n,n}(X)$ by bidegree. Thus $K^t=-K$.
The same argument applies to $J$ and gives $J^t=-J$.

Taking transposes in $J=[L_h,\Pi]$ now gives
\[
 -J=J^t=-[L_h,\Pi^t],
 \qquad\text{hence}\qquad
 [L_h,\Pi^t-\Pi]=0.
\]
Set $S=\Pi^t-\Pi$. Since the Poincar\'e pairing pairs complementary
degrees, $\Pi^t$, like $\Pi$, has total degree $-2$.
Thus $S$ lowers degree by two and commutes with $L_h$.
We show that these two properties force $S=0$.

Let $x\in H^j(X,\mathbb C)$ be an even primitive class, with
$j\leq n$. If $j=0$, then $Sx=0$ by degree. If $j\geq2$,
primitivity and commutation give
\[
 L_h^{n-j+1}Sx=S\bigl(L_h^{n-j+1}x\bigr)=0.
\]
Multiplying this equality once more by $L_h$ gives
$L_h^{n-j+2}Sx=0$. Hard Lefschetz makes this power injective on
$H^{j-2}(X,\mathbb C)$, so $Sx=0$. Since every even cohomology class
is a sum of Lefschetz powers of primitive classes, commutation
with $L_h$ implies that $S$ vanishes on all of $V$.
Therefore $\Pi^t=\Pi$, completing the proof.
\end{proof}

\begin{rem}\label{rem-operators}
The proof of Lemma~\ref{lem:contraction} gives, on the full Hodge
cohomology $\mathcal V$ and for every $z\in\mathbb C$,
\begin{equation}\label{eq:finite-conjugation}
 \begin{aligned}
 e^{-zL_h}De^{zL_h}
   &=D-z[L_h,D]+\frac{z^2}{2}[L_h,[L_h,D]]\\
   &=A_{a+z\eta+z^2c}
       +\bigl(K-z[L_h,\Pi]\bigr)+\Pi.
 \end{aligned}
\end{equation}
The three terms in the second line have bidegrees $(0,2)$,
$(-1,1)$ and $(-2,0)$, respectively. Thus the multiplication
component preserves the first Hodge index, while the other two
components lower it by one and two. We will use this decomposition
for the rank comparison in Section~\ref{sec:volume}.
\end{rem}

We now examine the consequences of the commutation relation
$[D,\overline D]=0$ for the three classes $a,\eta,c$.
As in degree one, this relation will yield orthogonality
with respect to the real part of a Hermitian form. Consider
\[
 \langle u,w\rangle_h
    =\int_X h^{n-2}\overline u\,w,
 \quad \textrm{where} \  u,w\in H^{0,2}(X).
\]
Every class in $H^{0,2}(X)$ is primitive, since its product with
$h^{n-1}$ vanishes by Hodge type. The Hodge-Riemann bilinear
relations therefore make this form positive definite.
We write $\|u\|_h^2=\langle u,u\rangle_h$.

The following lemma derives the required orthogonality from
the commutation relation and the contraction identities.
The resulting norm identity will imply that $D(e^{ih})=0$ forces
$\eta=K(h)=0$ and $a=c$.

\begin{lem}\label{lem:commuting-orthogonality}
Suppose that $D=A_a+K+\Pi$ commutes with $\overline D$ on $V$.
Then
\begin{equation}\label{eq:commutation-relations}
 K(\overline a)=\overline K(a),\qquad
 \Pi(\overline\eta)=\overline\Pi(\eta),\qquad
 \Pi(\overline a)=\overline\Pi(a).
\end{equation}
Consequently the three Hermitian products
$\langle a,\eta\rangle_h$, $\langle\eta,c\rangle_h$ and
$\langle a,c\rangle_h$ are real. In particular, $a-c$ and $i\eta$
are orthogonal for the real part of $\langle-,-\rangle_h$, and
\begin{equation}\label{eq:degree-two-pythagoras}
 \bigl\|e^{-ih}\cup D(e^{ih})\bigr\|_h^2
       =\|a-c\|_h^2+\|\eta\|_h^2.
\end{equation}
\end{lem}

\begin{proof}
We separate the commutator according to total cohomological degree.
The operators $A_a,K,\Pi$ have degree shifts $2,0,-2$, respectively,
and their conjugates have the same shifts. The components of
$[D,\overline D]=0$ with shifts $2,-2,0$ therefore give, respectively,
\[
 \begin{aligned}
 \relax[A_a,\overline K]+[K,A_{\overline a}]&=0,\\
 [K,\overline\Pi]+[\Pi,\overline K]&=0,\\
 [A_a,\overline\Pi]+[K,\overline K]
                         +[\Pi,A_{\overline a}]&=0.
 \end{aligned}
\]
Evaluate the first and third identities at $1_X$ and the second at
$h$. Contractions vanish on $1_X$, while $\Pi(h)=\overline\Pi(h)=0$
by bidegree. These evaluations give precisely
\eqref{eq:commutation-relations}. Each equality identifies its
left-hand side with its complex conjugate. Hence $K(\overline a)$
is a real class in $H^{1,1}(X)$, and $\Pi(\overline\eta)$ and
$\Pi(\overline a)$ are real scalars.

We now convert these reality statements into orthogonality.
The product $\overline a\,h^{n-1}$ is zero by Hodge type. Applying
the degree-zero derivation $K$ gives
\[
 0=K(\overline a)h^{n-1}
                  +(n-1)\overline a\,\eta h^{n-2}.
\]
After integration, this reads
\[
 (n-1)\langle a,\eta\rangle_h
             =-\int_X h^{n-1}K(\overline a)\in\mathbb R.
\]
For the other two products, the computation in
Lemma~\ref{lem:contraction} gives
$\Pi(h^n)=n(n-1)h^{n-2}c$. Using also $\Pi^t=\Pi$, we obtain,
for $u=a$ or $u=\eta$,
\[
 \begin{aligned}
 n(n-1)\langle u,c\rangle_h
   &=\int_X\overline u\,\Pi(h^n)\\
   &=\int_X\Pi(\overline u)h^n
     =\Pi(\overline u)\int_Xh^n\in\mathbb R.
 \end{aligned}
\]
In particular, $\langle a-c,\eta\rangle_h$ is real. Multiplication
of the second entry by $i$ makes this product purely imaginary, so
\[
 \operatorname{Re}\langle a-c,i\eta\rangle_h=0.
\]
Finally, Lemma~\ref{lem:contraction} gives
$e^{-ih}\cup D(e^{ih})=a+i\eta-c$. Expanding its squared norm
and using this orthogonality proves
\eqref{eq:degree-two-pythagoras}.
\end{proof}

Under the assumptions $[D,\overline D]=0$ and $D(e^{ih})=0$,
Lemma~\ref{lem:commuting-orthogonality} gives $K(h)=0$ and $a=c$.
These relations do not yet imply that $K$ and $\Pi$ vanish.
We therefore need positive quadratic expressions whose vanishing forces
the operators themselves to vanish. Poincar\'e duality provides a way
to construct them: it represents endomorphisms of $V$ by
cohomology classes on $X\times X$, where we can apply the
Hodge-Riemann bilinear relations. The crucial point will be
that the commutator identities of Lemma~\ref{lem:contraction}
translate into primitivity conditions for these classes. \\

Equip $X\times X$ with the ample class
\[
 \ell=c_1(L\boxtimes L)=p_1^*h+p_2^*h.
\]
Identify $\operatorname{End}(V)$ with the even-even K\"unneth
summand $V\otimes V$ by sending $u\otimes w$ to the operator
$x\mapsto P(w,x)u$, and denote by $C_T$ the cohomological
correspondence associated with $T$. If $T$ has bidegree $(r,s)$,
then $C_T$ has Hodge type $(n+r,n+s)$, because the Poincar\'e
dual of a class of type $(p,q)$ has type $(n-p,n-q)$.

Under this identification, composition with $L_h$ on the left
corresponds to multiplication by $p_1^*h$, while composition
on the right corresponds to multiplication by $p_2^*h$.
Consequently
\[
 C_{[L_h,T]}=(p_1^*h-p_2^*h)\cup C_T.
\]
The difference $p_1^*h-p_2^*h$ restricts to zero on the diagonal
and is not ample. To obtain primitivity with respect to the
ample sum $\ell$, we modify the correspondence by inserting a sign on
the second factor. Let $\epsilon$ act by $(-1)^j$ on
$H^{2j}(X,\mathbb C)$ and set
\[
                    \Gamma_T=(1\otimes\epsilon)C_T.
\]
Since $\epsilon L_h=-L_h\epsilon$, the commutator now corresponds
to multiplication by $\ell$. More precisely,
\begin{equation}\label{eq:product-kernels}
 \Gamma_{[L_h,T]}=\ell\Gamma_T,
 \qquad
 \int_{X\times X}\Gamma_T\overline{\Gamma_U}
      =(-1)^n\operatorname{Tr}\bigl(T(\overline U)^t\bigr).
\end{equation}
For the second identity, the correspondence identification (see
\cite[Proposition~1.3.6(ii)(c)]{Kle68}) gives
\[
 \int_{X\times X}C_T\overline{C_U}
      =\operatorname{Tr}\bigl(T(\overline U)^t\bigr).
\]
Since the correspondences lie in the even-even K\"unneth summand,
the supertrace in the usual formula is an ordinary trace.
Applying $\epsilon$ to both classes on the second copy of $X$
contributes $(-1)^n$
whenever their product has top degree there, proving
\eqref{eq:product-kernels}. 

We can now see why the Hodge-Riemann relations apply.
When $K(h)=0$, Lemma~\ref{lem:contraction} gives
$[L_h,K]=0$, and hence $\ell\Gamma_K=0$.
The same lemma gives
$\operatorname{ad}_{L_h}^3(\Pi)=0$, hence
$\ell^3\Gamma_\Pi=0$.
The classes $\Gamma_K$ and $\Gamma_\Pi$ have types
$(n-1,n+1)$ and $(n-2,n)$, respectively. Their degrees are
therefore $2n$ and $2n-2$, so these are precisely the
primitivity conditions on the $2n$-dimensional variety
$X\times X$.\\

The next proposition translates the resulting strict
Hodge-Riemann positivity into trace inequalities.
In Section~\ref{subsec:evaluation}, these nonnegative trace
expressions will add up to $\operatorname{Tr}(D\overline D)$
once $a=c$ and $K(h)=0$. Commutation and nilpotence will force
this trace to vanish, and strict positivity will then give
$K=\Pi=0$.

\begin{prop}\label{prop:product-positivity}
The bivector term and, when $K(h)=0$, the Kodaira-Spencer term satisfy
\[
\begin{aligned}
 2\operatorname{Tr}(A_c\overline\Pi)
    &=(-1)^n\int_{X \times X} \ell^2\Gamma_\Pi\,\overline{\Gamma_\Pi}\geq 0,\\
 \operatorname{Tr}(K\overline K)
    &=(-1)^{n+1}\int_{X \times X} \Gamma_K\,\overline{\Gamma_K}\geq 0
          \quad \text{if } K(h)=0.
\end{aligned}
\]
Equality holds in the first line precisely when $\Pi=0$ on $V$, and in the second precisely when $K=0$ on $V$.
\end{prop}
\begin{proof}
First, $\Gamma_\Pi$ has type $(n-2,n)$ on $X \times X$. By \eqref{eq:operators-L} and \eqref{eq:product-kernels},
\[
 \ell^2\Gamma_\Pi=2\Gamma_{A_c},\qquad \ell^3\Gamma_\Pi=0.
\]
Since $\deg\Gamma_\Pi=2n-2$, the second equation is the primitivity condition on the $2n$-fold $X \times X$. The Hodge-Riemann bilinear relations on $X \times X$ give the first inequality, with equality only when $\Gamma_\Pi=0$. Substituting the first equation in the integral and using $\Pi^t=\Pi$ in \eqref{eq:product-kernels} identifies the integral, with its sign, with $2\operatorname{Tr}(A_c\overline\Pi)$.\\

Next, $\Gamma_K$ has type $(n-1,n+1)$, so it is in the middle degree. If $K(h)=0$, then $[L_h,K]=0$ by \eqref{eq:operators-L}, and hence $\ell\Gamma_K=0$ by \eqref{eq:product-kernels}. Thus $\Gamma_K$ is primitive. Again, the Hodge-Riemann bilinear relations on $X \times X$ give the second inequality. In \eqref{eq:product-kernels}, the sign $(-1)^n$ together with $K^t=-K$ gives the stated identity with $\operatorname{Tr}(K\overline K)$. Equality holds precisely when $\Gamma_K=0$. The injectivity of $T\mapsto\Gamma_T$ gives both asserted vanishing criteria on $V$.
\end{proof}

\subsection{Injectivity in degree two}\label{subsec:evaluation}
To prove the injectivity of $\Theta^2_{F,i}$, we will show that if
the transported operator $D_\beta$ vanishes on the single class
$e^{ih}$, then $\beta=0$. As in degree one, we begin by using
positivity on $X$ to compare the terms in this evaluation. The new
step is to view the remaining components of $D_\beta$ as
correspondences on $X\times X$. There, positivity applies to the
operators themselves and forces them to vanish. \\

We start with the key vanishing lemma. Since Multiplication by a class in $H^{0,2}(Y)$ commutes with multiplication
by its complex conjugate, the cohomological Fourier-Mukai
transform is real, the transported operator on $X$ has the same
property. We combine the orthogonality of
Lemma~\ref{lem:commuting-orthogonality} with the positivity of
Proposition~\ref{prop:product-positivity} to prove the following:

\begin{lem}\label{lem:commuting-evaluation}
Let $h=c_1(L)$ for an ample line bundle $L$ on $X$. If $D\in\mathscr E_X$ satisfies $D\overline D=\overline D D$ and $D(e^{ih})=0$, then $D=0$ on $H^{\mathrm{ev}}(X,\mathbb C)$.
\end{lem}
\begin{proof}
If $\dim X\leq1$, every degree-two Hochschild operator acts as zero
on even cohomology by bidegree. Assume $n=\dim X\geq2$, write
$D=A_a+K+\Pi$, and put $\eta=K(h)$ and $c=\Pi(h^2)/2$.
The commutation hypothesis allows us to apply
Lemma~\ref{lem:commuting-orthogonality}. Since $D(e^{ih})=0$,
its norm identity gives
\[
                   0=\|a-c\|_h^2+\|\eta\|_h^2.
\]
Both terms are nonnegative, and the Hermitian form is positive
definite. Thus $\eta=0$ and $a=c$.

It remains to show that $K$ and $\Pi$ vanish. Each summand of $D$ lowers $p-q$ by two. Since $-n\leq p-q\leq n$, we have $D^{n+1}=0$, and commutation gives
\[
 (D\overline D)^{n+1}=D^{n+1}\overline D^{\,n+1}=0.
\]
Thus $\operatorname{Tr}(D\overline D)=0$. We compute this trace with
respect to total cohomological degree. Since $a=c$, the shifts of the
three components $A_c,K,\Pi$ are $2,0,-2$, and complex conjugation
preserves these shifts. Only products of total shift zero contribute,
so
\[
 \operatorname{Tr}(D\overline D)
   =\operatorname{Tr}(K\overline K)
      +\operatorname{Tr}(A_c\overline\Pi)
      +\operatorname{Tr}(\Pi A_{\overline c}).
\]
By cyclicity of the trace, the last two terms are complex conjugates.
Proposition~\ref{prop:product-positivity} makes
$\operatorname{Tr}(A_c\overline\Pi)$ real, so these terms are equal.
Consequently
\[
 0=\operatorname{Tr}(D\overline D)
   =\operatorname{Tr}(K\overline K)+2\operatorname{Tr}(A_c\overline\Pi).
\]
As $K(h)=0$, both summands are nonnegative by the same proposition. The characterizations of equality give $K=\Pi=0$ on even cohomology. Finally $a=c=\Pi(h^2)/2=0$, so all three components of $D$ vanish.
\end{proof}

\begin{proof}[Proof of Theorem~\ref{thm:twoforms}]
Proposition~\ref{prop:degree-one} treats degree one. Let
$\beta\in H^{0,2}(Y)$ lie in the kernel of $\Theta^2_{F,i}$. By
\eqref{eq:evaluation-map} at $z=i$, we have $D_\beta(e^{ih})=0$, since multiplication
by $e^{ih}$ is invertible. By \eqref{eq:transport}, $D_\beta|_V$ belongs
to $\mathscr E_X$. Moreover, $F$ is real and the even classes $\beta$
and $\bar\beta$ commute under cup product. Thus
\[
 [D_\beta,\overline{D_\beta}]=F^{-1}[m_\beta,m_{\bar\beta}]F=0.
\]
Lemma~\ref{lem:commuting-evaluation} gives $D_\beta=0$ on $V$. Since
$F$ identifies the even cohomology of $X$ with that of $Y$, it follows
that $m_\beta$ is zero on even cohomology. Evaluating at $1_Y$ gives
$\beta=0$. Applying injectivity to the inverse equivalence with an ample
line bundle on $Y$ gives the reverse dimension inequality. Hence
$\Theta^2_{F,i}$ is an isomorphism, and Hodge symmetry gives
$h^{2,0}(X)=h^{2,0}(Y)$.
\end{proof}

\begin{rem}\label{rem:degree-three-polarization}
In degree three, varying the parameter $z$ need not suffice to
obtain injectivity for a fixed polarization. Let
\[
                  X\subset\mathbb P^2\times\mathbb P^2
\]
be a general smooth hypersurface of bidegree $(3,3)$, and denote
by $H_1,H_2$ the restrictions of the two hyperplane classes.
This Calabi-Yau threefold contains rational curves $C_1,C_2$
of bidegrees $(1,0)$ and $(0,1)$, respectively
\cite[Proposition~7.4.2]{BvS95}. For a general $X$, the schemes
of lines of either bidegree are reduced and zero-dimensional.
Indeed, the defining equations generate the rank-four bundle
$\operatorname{Sym}^3(S^\vee)\boxtimes\mathcal O_{\mathbb P^2}(3)$
on $(\mathbb P^2)^\vee\times\mathbb P^2$, where $S$ is the
tautological rank-two bundle, so a general section vanishes
transversely; the other bidegree is treated by exchanging the
two factors. Thus $H^0(C_j,N_{C_j/X})=0$. Adjunction and the
splitting of vector bundles on $\mathbb P^1$ give
$N_{C_j/X}\simeq\mathcal O_{\mathbb P^1}(-1)^{\oplus2}$.
Thus $\mathcal O_X,\mathcal O_{C_1},\mathcal O_{C_2}$ are
spherical objects. Write $T_E$ for the spherical twist associated
with $E$ \cite{ST01}, and consider
\[
 \mathcal G=
 (-\otimes\mathcal O_X(1,0))
 \circ T_{\mathcal O_X}^{-1}
 \circ T_{\mathcal O_{C_2}}^{-1}
 \circ T_{\mathcal O_{C_1}}
 \circ T_{\mathcal O_X}.
\]
Let $\Phi=\mathcal G^{-1}$ and let $F$ be its cohomological
transform.

Put $p=[\mathrm{pt}]$ and $\delta=[C_2]-[C_1]$.
The Mukai vectors of the curve sheaves are
$v(\mathcal O_{C_j})=[C_j]+p$. Using the standard formula
\[
             [T_E(A)]=[A]-\chi(E,A)[E],
\]
one checks that the four twists in $\mathcal G$ send $p$
to $p+\delta$. Since $H_1\cup\delta=-p$, the final tensor
operation sends $p+\delta$ to $\delta$. Consequently
\[
                         F^{-1}(p)=\delta.
\]
Preservation of the Mukai pairing therefore gives, for any
ample class $h=c_1(L)$,
\[
 \bigl(F(e^{zh})\bigr)_0
       =\int_X e^{-zh}\delta
       =-z\int_X h\delta,
\]
where the subscript $0$ denotes the degree-zero component.

Moreover, $F$ acts as the identity on $H^{0,3}(X)$:
spherical twists act trivially on odd cohomology, and tensoring
with a line bundle fixes $H^{0,3}(X)$ by Hodge type.
Since $F(e^{zh})$ has diagonal Hodge type, only its degree-zero
component can have a nonzero product with
$\beta\in H^{0,3}(X)$. The definition of $\Theta^\bullet_{F,z}$
thus gives
\[
                 \Theta^3_{F,z}(\beta)
                    =-z\left(\int_X h\delta\right)\beta.
\]

For the ample polarization $h=H_1+H_2$, we have
$\int_X h\delta=0$. Hence $\Theta^3_{F,z}$ is the zero map
for every $z\in\mathbb C$, although $h^{0,3}(X)=1$.
If instead we choose the ample polarization $h'=2H_1+H_2$,
then $\int_X h'\delta=-1$, and the resulting map is
\[
                         \Theta^3_{F,z}=z\,\mathrm{id}.
\]
It is therefore injective for every $z\ne0$.

Thus, any attempt to obtain injectivity in degree three
by choosing a suitable complex parameter must also allow
for a suitable choice of polarization. This contrasts with
Theorem~\ref{thm:twoforms}, where $\Theta^2_{F,i}$ is injective
for every ample line bundle $L$.
\end{rem}

\section{Derived equivalences of stably holomorphically symplectic varieties}\label{sec:volume}

Theorem~\ref{thm:twoforms} compares the spaces of holomorphic
two-forms on derived-equivalent varieties. When one of these forms
is symplectic, its powers carry information in every degree.
Its top power is a volume form, and the ranks of cup product by
its successive powers determine all the Hodge numbers. We use these
observations to show that every derived partner of a stably
holomorphically symplectic variety is again of this kind and has
the same Hodge numbers. In odd dimension, the argument takes place
on the product with an elliptic curve.

We then ask whether this equality of Hodge numbers can be realized
by a single algebraic correspondence that also preserves cup product.
We construct such a correspondence from the Fourier-Mukai transform
and show that it induces an isomorphism of rational homological
motives as unital algebra objects. Under the composition-nilpotence
conjecture, the underlying isomorphism of motives lifts to rational
Chow motives.

\subsection{Equality of Hodge numbers}\label{subsec:hodge-numbers}

\begin{defn}\label{def:stably-symplectic}
A smooth projective complex variety $X$ is \emph{stably
holomorphically symplectic} if there exists a smooth projective
complex variety $Z$ such that $X\times Z$ admits a holomorphic
symplectic form.
\end{defn}

The following characterization also shows that the auxiliary
factor can always be chosen to be a point or an elliptic curve,
according to the parity of $\dim X$. We regard a point as
holomorphic symplectic.

\begin{prop}\label{prop:stable-characterization}
Let $X$ be a smooth projective complex variety of dimension $n$.
The following conditions are equivalent:
\begin{enumerate}
\item $X$ is stably holomorphically symplectic.
\item $K_X\simeq\mathcal O_X$, and the complex conjugate of a
holomorphic volume form belongs to the unital graded subalgebra of
$H^\bullet(X,\mathcal O_X)$ generated by
$H^1(X,\mathcal O_X)$ and $H^2(X,\mathcal O_X)$.
\item If $n$ is even, $X$ is holomorphic symplectic; if $n$ is odd,
$X\times E$ is holomorphic symplectic for every elliptic curve $E$.
\end{enumerate}
\end{prop}
\begin{proof}
Suppose first that $X\times Z$ admits a holomorphic symplectic
form $\Omega$. Restricting the trivial canonical bundle
$K_{X\times Z}=p_X^*K_X\otimes p_Z^*K_Z$ to a fibre of each
projection shows that $K_X$ and $K_Z$ are trivial. By K\"unneth,
\[
 \Omega=p_X^*\omega_X
       +\sum_j p_X^*\alpha_j\wedge p_Z^*\beta_j
       +p_Z^*\omega_Z,
\]
where $\omega_X,\omega_Z$ are holomorphic two-forms and
$\alpha_j,\beta_j$ are holomorphic one-forms on the respective
factors. Put $n+\dim Z=2N$ and choose volume forms $\tau_X,\tau_Z$.
Since $X\times Z$ is connected and projective,
\[
          \Omega^N=c\,p_X^*\tau_X\wedge p_Z^*\tau_Z,
          \qquad c\in\mathbb C^*.
\]
Every factor from $X$ in the expansion on the left is a product
of $\omega_X$ and the $\alpha_j$. Taking the coefficient of
$\tau_Z$ in the top K\"unneth component therefore expresses
$c\tau_X$ as a linear combination of such products. Complex
conjugation proves the second condition.

Assume the second condition and choose such a volume form $\tau_X$.
If $n=2d>0$, choose a nonzero monomial
of degree $2d$ in an expression for the conjugate volume form.
It contains an even number of degree-one factors. Grouping them
in pairs gives, up to sign, a product
$u_1\cdots u_d\ne0$ with $u_j\in H^2(X,\mathcal O_X)$.
The coefficient of $t_1\cdots t_d$ in
$(\sum_j t_ju_j)^d$ is $d!u_1\cdots u_d$, so there exists
$u\in H^2(X,\mathcal O_X)$ with $u^d\ne0$. The top form
$\bar u^d$ is a nonzero section of the trivial line bundle $K_X$,
hence vanishes nowhere. Thus $\bar u$ is holomorphic symplectic.

If $n$ is odd, let $E$ be an elliptic curve and let
$0\ne e\in H^1(E,\mathcal O_E)$. The product has trivial canonical
bundle, and the conjugate volume form $\bar\tau_X\boxtimes e$
is generated in degrees one and two. Applying the even-dimensional
argument to $X\times E$ proves the third condition. The case
$n=0$ is immediate. Finally, the third condition implies the first
by taking $Z$ to be a point or an elliptic curve.
\end{proof}

\begin{prop}\label{prop:stable-derived}
Let $X$ and $Y$ be derived-equivalent smooth projective complex
varieties. If $X$ is stably holomorphically symplectic, then
$\dim Y=\dim X$ and $Y$ is stably holomorphically symplectic.
\end{prop}
\begin{proof}
Write $n=\dim X$. By Proposition~\ref{prop:stable-characterization},
$K_X$ is trivial, so the Serre functor of $D^b(X)$ is the shift
$[n]$. Transporting it to $Y$ and evaluating on $\mathcal O_Y$
gives $\mathcal O_Y[n]\simeq K_Y[\dim Y]$. Hence $\dim Y=n$
and $K_Y\simeq\mathcal O_Y$.

Fix an ample line bundle on $X$ and let $F$ be the cohomological
transform of the equivalence. Write $A_X$ and $A_Y$ for the
unital graded subalgebras of $H^\bullet(X,\mathcal O_X)$ and
$H^\bullet(Y,\mathcal O_Y)$ generated in degrees one and two.
The homomorphism $\Theta^\bullet_{F,i}$ is surjective in these
degrees by Proposition~\ref{prop:degree-one} and
Theorem~\ref{thm:twoforms}. It therefore maps $A_Y$ onto $A_X$.
By Proposition~\ref{prop:stable-characterization}, a generator
$\tau\in H^n(X,\mathcal O_X)$ belongs to $A_X$. Choose
$\beta\in A_Y^n$ with $\Theta^n_{F,i}(\beta)=\tau$.
Then $\beta\ne0$, and since $K_Y$ is trivial, it generates
$H^n(Y,\mathcal O_Y)$. The same characterization proves that $Y$
is stably holomorphically symplectic.
\end{proof}

We now compare Hodge numbers, beginning with holomorphic symplectic
varieties. Write $m_u$ for cup product by $u$ and index the
Hochschild diagonals by
\[
              M_k(Z)=\bigoplus_{q-p=k}H^{p,q}(Z).
\]
Thus $m_u$ sends $M_k(Z)$ into $M_{k+2}(Z)$ when $u$ has type $(0,2)$.

\begin{lem}\label{lem:symplecticstrings}
Let $Z$ be a smooth projective variety of dimension $2d$, and let
$\beta\in H^{0,2}(Z)$ be the complex conjugate of a holomorphic
symplectic form. For $k\in\mathbb Z$ and $r\geq0$, write
\[
 R_Z(k,r)=\operatorname{rk}\bigl(
 m_{\beta^r}:M_k(Z)\longrightarrow M_{k+2r}(Z)\bigr).
\]
For every $p\in\mathbb Z$ and $0\leq q\leq d$, one has
\begin{equation}\label{eq:symplectic-rank-formula}
 R_Z(q-p,d-q)
 =h^{p,q}(Z)+\sum_{j=1}^{q}
   \bigl(h^{p-j,q-j}(Z)+h^{p+j,q-j}(Z)\bigr),
\end{equation}
where Hodge numbers outside the Hodge range are understood to be zero.
In particular, the ranks $R_Z(k,r)$ are independent of $\beta$
and determine all the Hodge numbers of $Z$.
\end{lem}

\begin{proof}
Put $L=m_\beta$. For $0\leq b\leq d$, symplectic linear algebra
gives an isomorphism of bundles
\[
 \bar\beta^{\,d-b}\wedge-:\Omega_Z^b
       \xrightarrow{\ \sim\ }\Omega_Z^{2d-b}.
\]
Taking sheaf cohomology and then complex conjugates gives, for every $a$,
\[
 L^{d-b}:H^{a,b}(Z)\xrightarrow{\ \sim\ }H^{a,2d-b}(Z).
\]
Consequently $L^r:H^{a,b}\to H^{a,b+2r}$ is injective when
$b+r\leq d$ and surjective when $b+r\geq d$. For injectivity,
$L^{d-b}$ factors through $L^r$. For surjectivity, when the target
is nonzero, the Lefschetz isomorphism from $H^{a,2d-b-2r}$ onto
$H^{a,b+2r}$ factors through $L^r:H^{a,b}\to H^{a,b+2r}$.

Fix $p$ and $q$ as in the statement. The summand of $M_{q-p}(Z)$
with first Hodge index $p+j$ is $H^{p+j,q+j}(Z)$. On this summand
the map whose rank we wish to compute is
\[
 L^{d-q}:H^{p+j,q+j}(Z)\longrightarrow
                 H^{p+j,2d-q+j}(Z).
\]
For $j\leq0$ this map is injective and its rank is $h^{p+j,q+j}$.
For $j\geq0$ it is surjective and its rank is
$h^{p+j,2d-q+j}=h^{p+j,q-j}$ by symplectic duality.
Thus in both cases the rank is $h^{p+j,q-|j|}$. The different
summands map into distinct Hodge summands of the target, so their
ranks add. Separating the term $j=0$ and grouping the terms indexed by
$j$ and $-j$ gives \eqref{eq:symplectic-rank-formula}. \\

Every rank with $0\leq r\leq d$ occurs in this formula by taking
$q=d-r$ and $p=q-k$. For $r>d$ the rank is zero. This proves that the ranks are
independent of $\beta$. Conversely, $R_Z(-p,d)=h^{p,0}(Z)$.
For $1\leq q\leq d$, every term in the sum on the right of
\eqref{eq:symplectic-rank-formula} has second index less than $q$.
Induction on $q$ therefore recovers $h^{p,q}(Z)$ for $q\leq d$.
The equality $h^{p,q}(Z)=h^{p,2d-q}(Z)$ completes the proof.\\
\end{proof}

\begin{prop}\label{prop:symplectichodge}
Let $X$ be a smooth projective holomorphic symplectic variety. Every smooth projective derived partner $Y$ is holomorphic symplectic, and $h^{p,q}(Y)=h^{p,q}(X)$ for all $p,q$.
\end{prop}

\begin{proof}
Fix an equivalence from $X$ to $Y$ with cohomological transform $F$
and write $\dim X=2d$. By Propositions~\ref{prop:stable-derived}
and~\ref{prop:stable-characterization}, $Y$ is also holomorphic
symplectic of dimension $2d$. The case $d=0$ is immediate, so
assume $d\geq1$.

Let $\tau_X\in H^{0,2}(X)$ be the conjugate of a holomorphic symplectic
form, and fix $h=c_1(L)$ for an ample line bundle $L$ on $X$. By
Theorem~\ref{thm:twoforms}, we may choose $\beta\in H^{0,2}(Y)$ with
$\Theta^2_{F,i}(\beta)=\tau_X$. By multiplicativity of the graded
homomorphism constructed in Section~\ref{subsec:graded}, we obtain
\[
                \Theta^{2d}_{F,i}(\beta^d)=\tau_X^d\ne0.
\]
It follows that $\bar\beta^d$ is a nonzero section of the trivial
line bundle $K_Y$. Since $Y$ is connected and projective, this
section vanishes nowhere. Hence $\bar\beta$ is holomorphic symplectic. \\

For the comparison on full cohomology, put $G=F\circ m_{e^{ih}}$ and
$T=G^{-1}m_\beta G$. Write $F^{-1}m_\beta F=A_a+K+\Pi$. By
\eqref{eq:finite-conjugation} at $z=i$ and the equality
$T(1_X)=\Theta^2_{F,i}(\beta)=\tau_X$, we have
\[
                T=m_{\tau_X}+\bigl(K-i[m_h,\Pi]\bigr)+\Pi.
\]
The three terms have bidegrees $(0,2)$, $(-1,1)$ and $(-2,0)$,
respectively.

We can now compare the ranks $R_X(k,r)$ and $R_Y(k,r)$ of
Lemma~\ref{lem:symplecticstrings}. Fix $k\in\mathbb Z$ and $r\geq0$.
Both $F$ and multiplication by $e^{ih}$ preserve the difference of
the Hodge indices, so $G$ identifies $M_k(X)$ with $M_k(Y)$. Since
$GT^r=m_\beta^rG$, the rank of $T^r$ on $M_k(X)$ is exactly
$R_Y(k,r)$.

On the other hand, $T^r-m_{\tau_X}^r$ strictly lowers the first Hodge
index: every term in the expansion of $T^r$ other than $m_{\tau_X}^r$
contains a factor that lowers this index, and no factor raises it.
Thus, with the Hodge summands ordered by their first index, the map
$T^r:M_k(X)\to M_{k+2r}(X)$ is block triangular. Its diagonal blocks
are
\[
       m_{\tau_X}^r:H^{p,p+k}(X)\longrightarrow H^{p,p+k+2r}(X).
\]
The rank of a block triangular map is at least the sum of the ranks
of its diagonal blocks. We conclude that
\begin{equation}\label{eq:symplectic-rank-bound}
                         R_X(k,r)\leq R_Y(k,r).
\end{equation}

Finally, apply the same argument to the inverse equivalence, starting
with the symplectic form $\bar\beta$ on $Y$. The resulting symplectic
form on $X$ may differ from $\overline{\tau_X}$, but
Lemma~\ref{lem:symplecticstrings} shows that it gives the same ranks
$R_X(k,r)$. We therefore obtain the reverse inequality for every
$k,r$. The ranks on the two varieties are equal, and
\eqref{eq:symplectic-rank-formula} recovers their Hodge numbers
row by row.

\end{proof}

\begin{prop}\label{prop:volumehodge}
Let $X$ and $Y$ be derived-equivalent smooth projective complex
varieties. If $X$ is stably holomorphically symplectic, then
$h^{p,q}(X)=h^{p,q}(Y)$ for all $p,q$.
\end{prop}
\begin{proof}
By Proposition~\ref{prop:stable-derived}, $Y$ is stably
holomorphically symplectic and has the same dimension as $X$.
If this dimension is even, both varieties are holomorphic symplectic
by Proposition~\ref{prop:stable-characterization}, and
Proposition~\ref{prop:symplectichodge} applies.

If the dimension is odd, choose an elliptic curve $E$.
The products $X\times E$ and $Y\times E$ are holomorphic symplectic
by Proposition~\ref{prop:stable-characterization}. They are also
derived equivalent: the exterior product of a kernel for the given
equivalence with $\mathcal O_{\Delta_E}$ defines an equivalence
of the products, with inverse obtained from the inverse kernel.
Thus Proposition~\ref{prop:symplectichodge} gives equality of their
Hodge numbers. For
$P_Z(u,v)=\sum_{p,q}h^{p,q}(Z)u^pv^q$, the K\"unneth formula reads
\[
 P_{X\times E}=P_X(1+u)(1+v),\qquad
 P_{Y\times E}=P_Y(1+u)(1+v).
\]
Cancelling the nonzero polynomial $(1+u)(1+v)$ proves $P_X=P_Y$.
\end{proof}

\begin{rem}\label{rem:hodge-structures}
The cohomological Fourier-Mukai transform need not preserve
cohomological degree, and the rank comparisons used above recover
the Hodge numbers without constructing degreewise rational
isomorphisms. In the next subsection, we normalize the transform
so that its degree-preserving component is an isomorphism induced
by an algebraic correspondence. After scalar normalization, this
isomorphism also preserves cup product. It gives an isomorphism of
homological motives as unital algebra objects and, consequently,
of rational Hodge structures in every degree.
\end{rem}

\subsection{Homological motives and cup product}\label{subsec:motives}
In this subsection we will prove the following:
\begin{volumerestated}
\volumestatement
\end{volumerestated}

We have already proved in Section~\ref{subsec:hodge-numbers} that $Y$
is stably holomorphically symplectic and has the same Hodge numbers
as $X$. It remains to construct the correspondence that preserves
cup product.
In even dimension, we first normalize the cohomological transform by
exponentials of divisor classes so that it preserves the line spanned
by the unit. The symplectic form and the Hochschild action then show
that this transform and its inverse do not increase degree. Their
degree-preserving components are therefore inverse and are induced by
the codimension-$\dim X$ parts of the normalized correspondences.
After rescaling to preserve the unit, the resulting isomorphism
intertwines cup product by each class with cup product by its image,
as we show using the full Hochschild action. In odd dimension, the
construction descends from the product with an elliptic curve.

We write $[-]_j$ for the component of cohomological degree $j$ and
$b_2(Z)=\dim H^2(Z,\mathbb Q)$.

\begin{lem}\label{lem:exponential-normalization}
Let $X,Y$ be smooth projective holomorphic symplectic varieties of
the same positive dimension, with $b_2(X)=b_2(Y)$, and let $F$ be
the cohomological transform of a derived equivalence. If $a$ is a
rational divisor class on $X$ and $r=[F(e^a)]_0\ne0$, then
\[
 F(e^a)=r e^b,
 \qquad \textrm{with} \  b=\frac{[F(e^a)]_2}{r}
             \in H^2(Y,\mathbb Q)\cap H^{1,1}(Y).
\]
\end{lem}
\begin{proof}
Let $n=\dim X$ and let $\mathfrak g_Z$ denote the rational LLV
algebra of $Z$. We use only that it contains the grading operator
$\mathsf H_Z|_{H^j(Z)}=(j-n)\operatorname{id}$ and all
multiplications by classes in $H^2(Z,\mathbb Q)$, and that
\begin{equation}\label{eq:llv-unit-image}
       \mathfrak g_Z\cdot1_Z=H^0(Z,\mathbb Q)\oplus H^2(Z,\mathbb Q).
\end{equation}
The last equality follows from the three-grading of $\mathfrak g_Z$
and the identification of its degree-two part with $H^2(Z,\mathbb Q)$
\cite[Proposition~3.45 and the proof of Corollary~3.53]{Lei26};
here $n>0$ ensures that $\mathsf H_Z(1_Z)=-n1_Z$ spans $H^0(Z)$.
Taelman's theorem gives $F\mathfrak g_XF^{-1}=\mathfrak g_Y$
\cite[Theorem~A]{Tae23}.

Put $v=F(e^a)$, so that $v_0=r1_Y$. If $v=r e^b$, comparison
in degree two forces $b=v_2/r$. We therefore choose this value of
$b$, which has the asserted rationality and Hodge type since $F$
is rational and preserves the difference of the Hodge indices.
Set $B=m_{e^{-b}}Fm_{e^a}$ and $t=B(1_X)=e^{-b}\cup v$.
Expanding $e^{-b}=1-b+\cdots$ in degrees zero and two gives
\[
                    t_0=r1_Y,\qquad t_2=v_2-rb=0.
\]
Thus the formula for $b$ cancels precisely the degree-two component;
we must still show that all higher-degree components of $t$ vanish.
Conjugation by $m_{e^a}$ and $m_{e^{-b}}$
preserves the corresponding LLV algebras, because $m_a,m_b$
belong to them. Consequently
\[
 W:=\mathfrak g_Y\cdot t
    =B(\mathfrak g_X\cdot1_X)
    =B\bigl(H^0(X,\mathbb Q)\oplus H^2(X,\mathbb Q)\bigr)
\]
has dimension $1+b_2(Y)$ and contains $t$.

The projection $W\to H^0(Y,\mathbb Q)\oplus H^2(Y,\mathbb Q)$
is surjective: $\mathsf H_Yt$ projects to $-nr1_Y$, while
$m_ct$ projects to $rc$ for every $c\in H^2(Y,\mathbb Q)$.
It is therefore an isomorphism. The vector
$(\mathsf H_Y+n\operatorname{id})t$ belongs to $W$. Its degree-zero
component is zero, and its degree-two component is $2t_2=0$ by
our choice of $b$. Its projection is therefore zero, so the vector
itself vanishes. Since $\mathsf H_Y+n\operatorname{id}$
acts by $j$ in degree $j$, this forces $t=r1_Y$, as required.
\end{proof}

For the remaining steps, we use the spaces of Hochschild operators
\[
 \mathscr C_Z^k=\{\iota_\gamma:\gamma\in HT^k(Z)\}
       \subset\operatorname{End}H^\bullet(Z,\mathbb C).
\]
By \eqref{eq:transport}, the cohomological transform $F$ conjugates
$\mathscr C_X^k$ to $\mathscr C_Y^k$. Conjugation by $m_{e^c}$
preserves $\mathscr C_Z^k$ for every rational divisor class $c$:
for integral multiples of a first Chern class, this is the action
of tensoring by a line bundle, and the general case follows by
polynomiality, as in Proposition~\ref{prop:graded-homomorphism}.
Thus every normalization $G=m_{e^{-b}}Fm_{e^a}$ by rational divisor
classes satisfies
\begin{equation}\label{eq:normalized-hochschild-transport}
                     G\mathscr C_X^kG^{-1}=\mathscr C_Y^k.
\end{equation}

\begin{lem}\label{lem:unit-filtration}
Let $X,Y$ be smooth projective holomorphic symplectic varieties of
dimension $2d>0$. Suppose that a rational linear isomorphism
$G:H^\bullet(X,\mathbb Q)\to H^\bullet(Y,\mathbb Q)$ preserves
the difference of the Hodge indices and satisfies, after
complexification,
\[
 G\mathscr C_X^2G^{-1}=\mathscr C_Y^2 \quad \textrm{and} 
 \quad G(1_X)=r1_Y, \ \text{with }r\ne0.
\]
Then $G$ and $G^{-1}$ do not increase cohomological degree.
Their degree-preserving components are inverse isomorphisms.
\end{lem}
\begin{proof}
Choose a holomorphic symplectic form $\sigma$ on $X$ and its inverse
bivector $\pi$. Symplectic linear algebra gives
\[
 [m_\sigma,\iota_\pi]=\mathsf H_{p,X},
 \qquad \mathsf H_{p,X}|_{H^{p,q}(X)}=(p-d)\operatorname{id},
 \qquad \iota_\pi(\sigma^d)=d\sigma^{d-1};
\]
see \cite[Lemma~2.9]{Tae23} for the commutator identity.
Since $m_{\bar\sigma}\in\mathscr C_X^2$ and $G$ commutes with
complex conjugation, we have
\[
 A:=Gm_\sigma G^{-1}=m_\eta+R+S,
\]
where $\eta\in H^{2,0}(Y)$ and $R,S$ have bidegrees $(1,-1)$
and $(0,-2)$. They vanish on $H^{\bullet,0}(Y)$, so induction
gives $G(\sigma^j)=r\eta^j$. In particular $\eta^d\ne0$.
As $K_Y$ is trivial, this top power vanishes nowhere and $\eta$
is symplectic.

On the other hand, $B:=G\iota_\pi G^{-1}$ belongs to
$\mathscr C_Y^2$ and kills $1_Y$. Hence $B=K'+\Pi'$, with
$K'$ a Kodaira-Spencer contraction and $\Pi'$ a bivector
contraction. The identity on the top power gives
\[
 B(\eta^d)=d\eta^{d-1}.
\]
Separating Hodge types yields $K'(\eta^d)=0$ and
$\Pi'(\eta^d)=d\eta^{d-1}$. Contraction with the volume form
$\eta^d$ identifies $H^1(Y,T_Y)$ with $H^{2d-1,1}(Y)$ and
$H^0(Y,\wedge^2T_Y)$ with $H^{2d-2,0}(Y)$. Thus $K'=0$ and
$\Pi'$ is contraction by the inverse bivector of $\eta$.
Consequently
\begin{equation}\label{eq:grading-lowering}
 G\mathsf H_{p,X}G^{-1}=[A,B]=\mathsf H_{p,Y}+N,
 \qquad N=[R+S,\Pi'].
\end{equation}
The two components of $N$ have bidegrees $(-1,-1)$ and $(-2,-2)$.

For $0\ne x\in H^{p,q}(X)$, let $p'$ be the largest first Hodge
index occurring in $Gx$. The component of that index in
$(\mathsf H_{p,Y}+N)Gx=(p-d)Gx$ gives $p'=p$.
Since $G$ preserves the difference of the Hodge indices, it cannot
increase total degree. The same argument applies to $G^{-1}$,
and their degree-preserving components are inverse.
\end{proof}

The normalized transform therefore has an invertible
degree-preserving component. To prove that its unital normalization
preserves cup product, we retain the transport of the full
Hochschild action on this component. Complex conjugation supplies
multiplication by a holomorphic symplectic form, and commutators
with contractions then recover multiplication by every cohomology
class. This is the content of the next proposition.

\begin{prop}\label{prop:multiplicative-normalization}
Let $F:H^\bullet(X,\mathbb Q)\to H^\bullet(Y,\mathbb Q)$ be the
cohomological transform of a derived equivalence between smooth
projective holomorphic symplectic varieties. Let $a,b$ be rational
divisor classes on $X,Y$, respectively, and put $G=m_{e^{-b}}Fm_{e^a}$.
Suppose that $G$ and $G^{-1}$ do not increase cohomological degree
and that $G(1_X)=r1_Y$ with $r\ne0$. If $G_0$ is the component
preserving degree, then $U=r^{-1}G_0$ is an isomorphism of unital
graded $\mathbb Q$-algebras, preserving the Hodge structures.
\end{prop}
\begin{proof}
Since $G$ and $G^{-1}$ do not increase degree, $U=r^{-1}G_0$ is
a graded rational linear isomorphism with $U(1_X)=1_Y$. It preserves
each Hodge bidegree, because $G$ preserves the difference of the
Hodge indices. We will prove that
\[
                 Um_xU^{-1}=m_{U(x)}
                 \qquad\text{for every }x\in H^\bullet(X,\mathbb C).
\]

First, $U$ transports contractions of each bidegree. Indeed, for
$\gamma\in H^q(X,\wedge^pT_X)$, the component of total degree $q-p$
of $G\iota_\gamma G^{-1}$ is $U\iota_\gamma U^{-1}$, since $G$ and
$G^{-1}$ only lower degree. By
\eqref{eq:normalized-hochschild-transport}, this component belongs
to $\mathscr C_Y^{p+q}$, which is graded by total degree. Hence
\[
 U\iota_\gamma U^{-1}=\iota_{\gamma'}
 \qquad\text{for some }\gamma'\in H^q(Y,\wedge^pT_Y).
\]

Choose a holomorphic symplectic form $\sigma$ on $X$. Applying this
observation to $\bar\sigma\in H^2(X,\mathcal O_X)$ shows that
$Um_{\bar\sigma}U^{-1}$ is multiplication by $U(\bar\sigma)$,
as follows by evaluation at $1_Y$. Since $U$ is rational, complex
conjugation gives
\[
                 Um_\sigma U^{-1}=m_\eta,
                 \qquad\eta=U(\sigma)\in H^{2,0}(Y).
\]

Now fix $x\in H^{p,q}(X)$. The bundle isomorphism
$\sigma^\flat:\wedge^pT_X\to\Omega_X^p$, induced by
$v\mapsto\iota_v\sigma$, gives
$x=H^q(\sigma^\flat)(\gamma)$ for some
$\gamma\in H^q(X,\wedge^pT_X)$. We have
\begin{equation}\label{eq:iterated-symplectic-commutator}
 (\operatorname{ad}m_\sigma)^p(\iota_\gamma)
       =(-1)^{p(p+1)/2}p!\,m_x.
\end{equation}
To see this, write a local polyvector as
$v_1\wedge\cdots\wedge v_p$ and use
$[m_\sigma,\iota_v]=-m_{\iota_v\sigma}$. The iterated commutator
replaces each contraction once by exterior multiplication, giving
$p!$ terms; reversing the order of the resulting one-forms gives
the displayed sign. This identity of sheaf maps passes to cohomology
without an additional sign, since $m_\sigma$ has even degree.
The same identity holds for any holomorphic two-form, including $\eta$.

Conjugating \eqref{eq:iterated-symplectic-commutator} by $U$ and
using the same identity on $Y$ therefore gives
\[
                   Um_xU^{-1}=m_{x'},
                   \qquad x'=H^q(\eta^\flat)(\gamma').
\]
Evaluation at $1_Y$ yields $x'=U(x)$, proving the desired identity
for every Hodge component and hence for every class. Applying it to
$U(y)$ now gives
\[
          U(x\cup y)=(Um_xU^{-1})U(y)=U(x)\cup U(y),
\]
which proves the proposition.
\end{proof}

We work with rational homological motives defined using Betti
cohomology: correspondences are algebraic cycles with rational
coefficients, modulo the kernel of the Betti cycle class map.
The multiplication on $h_{\mathrm{hom}}(Z)_{\mathbb Q}$ is induced
by the small diagonal of $Z^3$.

\begin{proof}[Proof of Theorem~\ref{thm:volume}]
Fix an equivalence with cohomological transform $F$, and choose
$h=c_1(L)$ for an ample line bundle $L$ on $X$.
By Proposition~\ref{prop:stable-derived}, $Y$ is stably holomorphically
symplectic of the same dimension as $X$. Suppose first that
$n=\dim X=2d>0$. Proposition~\ref{prop:stable-characterization}
shows that both varieties are holomorphic symplectic, and
Proposition~\ref{prop:volumehodge} gives equality of their second
Betti numbers. As in Proposition~\ref{prop:symplectichodge},
choose $\beta\in H^{0,2}(Y)$ such that
\[
                     \Theta^n_{F,i}(\beta^d)\ne0.
\]
Consider the rational polynomial $r(z)=[F(e^{zh})]_0$.
Since $F(e^{zh})$ is of diagonal Hodge type,
\[
                  \beta^d\cup F(e^{zh})=r(z)\beta^d.
\]
The nonvanishing above implies $r(i)\ne0$. Choose an integer $m>0$
with $r=r(m)\ne0$, and put $a=mh$.

Lemma~\ref{lem:exponential-normalization} gives
$F(e^a)=r e^b$. The class
$b=[F(e^a)]_2/r$ is rational of type $(1,1)$, hence is a rational
divisor class. The isomorphism
\begin{equation}\label{eq:normalized-motivic-transform}
                       G=m_{e^{-b}}F m_{e^a}
\end{equation}
satisfies $G(1_X)=r1_Y$, preserves the difference of the Hodge
indices, and transports the Hochschild action by
\eqref{eq:normalized-hochschild-transport}.
By Lemma~\ref{lem:unit-filtration}, $G$ and $G^{-1}$ do not
increase cohomological degree. If $G_0$ is the degree-preserving
component of $G$, Proposition~\ref{prop:multiplicative-normalization}
shows that $U=r^{-1}G_0$ is a unital graded algebra isomorphism
preserving the Hodge structures.

Choose rational algebraic correspondences
$C\in\mathrm{CH}^\bullet(X\times Y)_{\mathbb Q}$ and
$D\in\mathrm{CH}^\bullet(Y\times X)_{\mathbb Q}$ inducing $G$ and
$G^{-1}$: these are obtained from the Mukai vectors of the two
Fourier-Mukai kernels by multiplication with the divisor
exponentials in \eqref{eq:normalized-motivic-transform}.
Their codimension-$n$ components induce $G_0$ and $G_0^{-1}$.
Thus
\[
                 \Gamma=r^{-1}C_n,\qquad \Gamma'=rD_n
\]
induce $U$ and $U^{-1}$.

By the K\"unneth formula and Poincar\'e duality, an algebraic
correspondence is determined modulo homological equivalence by
its action on full cohomology. Hence
$\Gamma'\circ\Gamma=\Delta_X$ and
$\Gamma\circ\Gamma'=\Delta_Y$ modulo homological equivalence.
The same observation on $X\times X\times Y$, together with
the multiplicativity of $U$, gives
\[
 \Gamma\circ\Delta_X^{\mathrm{sm}}
       =\Delta_Y^{\mathrm{sm}}\circ(\Gamma\otimes\Gamma)
       \quad\text{modulo homological equivalence}.
\]
Here $\Delta_Z^{\mathrm{sm}}$ denotes the correspondence inducing
cup product. Since $U(1_X)=1_Y$, the unit is preserved as well.
This proves the claimed isomorphism of algebra objects in even
dimension.

Now let $n$ be odd. Choose an elliptic curve $E$ and a class
$h_E=c_1(M)$ with $M$ ample. By Proposition~\ref{prop:stable-characterization}, the products
$X\times E$ and $Y\times E$ are holomorphic symplectic. They are
also derived equivalent, as in Proposition~\ref{prop:volumehodge}.
Apply the preceding construction to the product equivalence and
the ample class $\widetilde h=p_X^*h+p_E^*h_E$.
Its cohomological transform is $F_E=F\otimes\operatorname{id}$.
For the integer $m$ supplied by that construction, put
$v_m=F(e^{mh})$ and $r=[v_m]_0\ne0$. The normalizing class on
the product satisfies
\[
 \widetilde b=p_Y^*b+m p_E^*h_E,
 \qquad b=\frac{[v_m]_2}{r},
\]
because $F_E(e^{m\widetilde h})=v_m\otimes e^{mh_E}$.
The normalized transform therefore factors as
\[
 m_{e^{-\widetilde b}}F_E m_{e^{m\widetilde h}}
       =\underbrace{\bigl(m_{e^{-b}}F m_{e^{mh}}\bigr)}_{G}
                       \otimes\operatorname{id}_{H^\bullet(E,\mathbb Q)}.
\]
It and its inverse do not increase degree. Applying them to
$x\otimes1_E$ shows that $G$ and $G^{-1}$ have the same property.
Their degree-preserving components satisfy
$(G\otimes\operatorname{id})_0=G_0\otimes\operatorname{id}$.
Consequently the unital algebra isomorphism constructed on the
products is
\[
                  U_{\mathrm{prod}}
                     =(r^{-1}G_0)\otimes\operatorname{id}.
\]
Restricting to classes $x\otimes1_E$ proves that $U=r^{-1}G_0$
is multiplicative. Choose algebraic correspondences $C,D$ inducing
$G,G^{-1}$ as above. The correspondences $r^{-1}C_n,rD_n$ on
$X\times Y$ and $Y\times X$ therefore give the required
isomorphism of algebra objects exactly as above. In dimension
zero both varieties are points.

Finally, the degree-preserving algebraic correspondences inducing
$U$ and $U^{-1}$ preserve each Hodge type. They give an isomorphism
of graded rational cohomology algebras compatible with the Hodge
structures, and in particular the stated equality of Hodge numbers.
\end{proof}

\subsection{A conditional lift to Chow motives}\label{subsec:chow}

The proof of Theorem~\ref{thm:volume} gives two Chow correspondences
whose compositions agree with the diagonals modulo homological
equivalence. To obtain an isomorphism of the underlying Chow motives,
it suffices to correct the inverse correspondence. We use the
following form of the nilpotence conjecture for algebraic cycles
\cite[Section~3.3]{KS16}.

\begin{conj}[Composition nilpotence]\label{conj:composition-nilpotence}
Let $Z$ be a smooth projective complex variety of dimension $d$.
If $R\in\mathrm{CH}^d(Z\times Z)_{\mathbb Q}$ has zero Betti
cycle class, there exists an integer $N\geq 1$ such that
$R^{\circ N}=0$ in $\mathrm{CH}^d(Z\times Z)_{\mathbb Q}$, where
the power is taken under composition of correspondences.
\end{conj}

\begin{rem}
Each of three classical conjectures implies this form of nilpotence.
\begin{itemize}
\item A homologically trivial correspondence $R$ is numerically trivial.
Voevodsky's smash-nilpotence conjecture therefore gives
$R^{\boxtimes N}=0$ for some $N$. Pulling this external power back
along the intermediate diagonals and pushing forward to the two
outer factors gives $R^{\circ N}=0$
\cite[Lemma~2.7 and Conjecture~4.2]{Voe95}.

\item The Kimura-O'Sullivan conjecture predicts that every rational
Chow motive is finite-dimensional. Kimura's nilpotence theorem
then applies directly to $R$ as a numerically trivial endomorphism
of $h_{\mathrm{CH}}(Z)_{\mathbb Q}$
\cite[Proposition~7.5]{Kim05}. See also \cite{Jan07}.

\item Finally, the Bloch-Beilinson conjecture gives a finite descending
filtration $F^\bullet$ on $\mathrm{CH}^d(Z\times Z)_{\mathbb Q}$.
Its first step consists of the homologically trivial cycles, and
composition sends $F^a\times F^b$ into $F^{a+b}$
\cite[Conjecture~2.1 and p.~292]{Jan94}. Thus $R^{\circ N}\in F^N=0$
for $N$ sufficiently large.
\end{itemize}
\end{rem}

\begin{prop}\label{prop:chow-lift}
\chowliftstatement
\end{prop}

\begin{proof}
Let $\Gamma$ and $\Gamma'$ be the Chow correspondences constructed
in the proof of Theorem~\ref{thm:volume}. They induce mutually
inverse isomorphisms of homological motives. Hence
\[
 R_X=\Delta_X-\Gamma'\circ\Gamma,\qquad
 R_Y=\Delta_Y-\Gamma\circ\Gamma'
\]
have zero Betti cycle class. By Conjecture~\ref{conj:composition-nilpotence},
both are nilpotent under composition. Thus
$\Gamma'\circ\Gamma=\Delta_X-R_X$ and
$\Gamma\circ\Gamma'=\Delta_Y-R_Y$ are invertible, with inverses
given by finite geometric sums. It follows that
$(\Gamma'\circ\Gamma)^{-1}\circ\Gamma'$ is a left inverse
of $\Gamma$, and $\Gamma'\circ(\Gamma\circ\Gamma')^{-1}$
is a right inverse. The two coincide, so $\Gamma$ itself gives the
required lift.
\end{proof}

\bibliographystyle{alpha}
\bibliography{Derived-invariance-h20-bibli}
\end{document}